\documentclass[12pt]{article}
\usepackage{MnSymbol}
\usepackage{amsmath,amsfonts,
verbatim}
\usepackage[all,cmtip]{xy}

\newcommand{\A}{\mathbb{A}}
\newcommand{\Q}{{\mathbb Q}}
\newcommand{\C}{{\mathbb C}}
\newcommand{\N}{{\mathbb N}}
\newcommand{\Z}{{\mathbb Z}}

\newcommand{\PP}{{\mathbf P}}
\newcommand{\HH}{\mathbb{H}}

\newcommand{\F}{\mathrm{F}}
\newcommand{\U}{{\mathbb U}}

\newcommand{\X}{{\mathbb X}}
\newcommand{\Y}{\mathbb{Y}}

\newcommand{\gh}{\mathbf{g}}

\newcommand{\pr}{{\rm pr}}

\newcommand{\la}{\langle}
\newcommand{\ra}{\rangle}
\newtheorem{pkt}{}[section]  
\newcommand{\bpk}{\begin{pkt}\rm }  
\newcommand{\epk}{\end{pkt}} 
\newcommand{\inv}{^{-1}}   

\newcommand{\be}{\begin{equation}}  
\newcommand{\ee}{\end{equation}}

\newcommand{\trd}{\mathrm{tr.deg}}

\newcommand{\acl}{{\rm acl}}

\newcommand{\LL}{\mathcal{L}}

\newcommand{\kk}{\mathrm{k}}
\newcommand{\G}{\mathrm{G}}\newcommand{\GG}{\mathrm{G}}
\newcommand{\GL}{\mathrm{GL}}
\newcommand{\SL}{\mathrm{SL}}

\newcommand{\E}{\mathrm{E}}

\newcommand{\Aut}{\mathrm{Aut}}
\newcommand{\Gal}{\mathrm{Gal}}
\newcommand{\tp}{\mathrm{tp}}
\newcommand{\sss}{\mathbf{s}}
\newcommand{\ttt}{\mathbf{t}}
\newcommand{\rrr}{\mathbf{d}}

\newcommand{\cl}{\mathrm{cl}}
\usepackage{graphicx}

\usepackage{tikz}

\title{Modular curves  and their pseudo-analytic cover}

\author{C.Daw\footnote{Supported by EPSRC New Investigator award EP/S029613/1} 
and B.Zilber\footnote{Supported by EPSRC Program grant ``Symmetries and Correspondences"}     }

\begin{document}

\maketitle

\section{Introduction}
\bpk {\bf Motivation and origins.}
The original aim of the paper is the extension of the project {\em categoricity of non-elementary theories of analytic covers} from  ``abelian" cases, such as $\exp: \C\to \mathbb{G}_m$ or $\exp_\Lambda: \C\to \E_\Lambda \ (\mbox{elliptic curve},$ to hyberbolic curves and possibly wider, see e.g. the survey part in \cite{Zspecial} for some history and references.

The first obstruction for this project is that an adequate formalism (that is the language) in which such an analytic cover can be properly presented is not easy to determine. Some attempts in this direction were by M.Gavrilovich \cite{Misha}, as well as later attempts by A.Harris \cite{Adam}. The one which we found satisfactory and applied here is based on the formalism close to the one applied in the recent  \cite{AZ}.  

More specifically, we formalise e.g. the case   $\exp: \C\to \mathbb{G}_m$ as a structure with two sorts $\U$ and $\F,$ where $ \U$ is the complex numbers with the $\Q$-module structure and the distinguished subgroup $2\pi i\Z,$ and $\F$ is the complex numbers as the field. Then $\exp:  \U\to \mathbb{G}_m(\F)$ is the map between the two sorts. But in fact along with
$\exp$ we automatically get in abelian cases the family 
$$\exp_n:  \U\to \mathbb{G}_m(\F),\ \ \exp_n: x\mapsto \exp(\frac{x}{n})$$
which of course agrees with the system of finite covers
$$\pr_{k,n}: \mathbb{G}_m(\F)\to \mathbb{G}_m(\F); \ y\mapsto y^\frac{k}{n},\mbox{ when } n|k.$$
 This we describe axiomatically by simple $L_{\omega_1,\omega}$ sentences. In terms of model-theoretic classification such a structure is a {\em fusion between a locally modular structure $\U$ and an algebraically closed field $\F.$} The case $\exp_\Lambda: \C\to \E_\Lambda $ similarly represents a fusion between a locally modular structure $\U_\Lambda,$  and an algebraically closed field $\F,$ where $\U_\Lambda$ is a $\Q$-module with distinguished $\Z$-module $\Lambda$ with an alternating form on it (accounting for the Weil pairing).  
\epk
\bpk \label{(i-ii)}
It is not clear a priori why such $L_{\omega_1,\omega}$-theories should be categorical in uncountable cardinals and the fact that they are must be of some significance.

The geometric value of the project is perhaps in the fact that the formulation of the categorical theory of the universal cover of a variety $\X$ (essentially the description of $\U$)  is  a formulation of a complete formal {\bf invariant} of $\X.$ 
By its  nature such an $L_{\omega_1,\omega}$-invariant is of ``algebraic type'' and the fact that it is equivalent to a notion given in topological/analytic terms indicates a possibility of connection to certain key conjectures of algebraic geometry such as the Hodge conjecture. 

Indeed, the most interesting outcome of the earlier works was establishing an equivalence between categoricity of the cover of 
$\E_\Lambda$ and the conjuction of the  two arithmetic facts:

 (i) the complete classification of the Galois action on the torsion of 
 $\E_\Lambda$ (the open image theorem by J.-P.Serre) and 
 
 (ii) the Kummer theory for 
  $\E_\Lambda$ (M.Bashmakov and K.Ribet)  

For abelian varieties the success of the program depends on the  extension of analogues of (i) and (ii) to abelian varieties, and  (ii) is known due to K.Ribet and M.Larsen. However, an analogue of Serre's theorem for abelian varieties is an open problem
  and therefore the best categoricity result here is under the assumption that the language names points of {\em the kernel} of the exponential map. This autmomatically removes the problem of determining the Galois action on the torsion points at the cost of weakening the formal invariant of $\X/\kk,$ the theory of the cover of $\X/\kk,$ to the formal invariant of $\X/K$ where $K$ is the field obtained from $\kk$ by adjoining all torsion points of $\X.$

The case of hyperbolic curves $\X$ was first considered by  
A.Harris \cite{Adam} and C.Daw and A.Harris  \cite{DH} in the context of modular curves $\Y(\Gamma)$ in the formalism (the choice of the language) which has names for each element of $\GL_2^+(\Q)$ acting on the upper half plane $\HH$ by M\"obius transformations. The proof of categoricity in this setting required essentially 
Serre's open image theorem for products of non-CM elliptic curves. Interestingly, in the analysis of $(\HH,j_N, \Y(N))$ Serre's theorem plays rather the role of (ii), while (i) is not needed since naming elements of $\GL_2(\Q)$ fixes the special points (CM-points). 
 This has a cost: one works out the formal invariant of $\Y(\Gamma)$ over the extension of the natural field of definition by special points.

\epk

\bpk {\bf Our setting.} Our current interest is  the case of the universal cover of hyperbolic curves, e.g. 
 the complex curve
 $\X=\PP^1\setminus \{ 0,1,\infty\}.$ However, before approaching this case we set ourselves a simpler task  of the {\em cover of the modular curve $\Y(N)$   universal in the class of modular curves}, which means that our structure incorporates the analytic covering maps $$j_n: \HH\to \Y(n), \mbox{ for all $n$ such that } N|n$$ agreeing with the algebraic finite covers $$\pr_{n,m}: \Y(n)\to \Y(m),\mbox{ for all $n,m$ such that } N|m,\ m|n.$$ 

In fact it is enough to classify the case $N=1,$ $\Y(1)=\Y=\mathbf{A}^1$ the affine line (the fact that for some  $n$ the covers $j_n: \HH\to \Y(n)$ are ramified does not matter in our setting).  

Note that $\Y(2)=\PP^1\setminus \{ 0,1,\infty\}$ when we consider $\Y(2)$ as an algebraic curve (without the level structure).

The important difference with the case of the proper universal cover is that, instead of the profinite completion $\hat{\Gamma}(N)$ of the respective fundamental group,  in the modular setting one has the group $\tilde{\Gamma}(N),$  the completion in the topology based on congruence subgroups, which for $N=1$ gives us
$$\tilde{\Gamma}(1)=\tilde{\Gamma}= \mathrm{SL}(2,\hat{\Z}),$$
where $$\hat{\Z}=\lim \Z/n\Z,$$
the projective limit of residue rings. 

\epk

\bpk The key problem, similar to earlier cases, is in classifying the saturated version of $\HH$ in the structure $(\HH, j_n,\Y(n))_{n\in \N}$  which is essentially reducible to understanding the structure on the
 projective limit 
 $$\tilde{\HH}:= \lim_{\leftarrow} \Gamma(n)\backslash \HH\cong  \lim_{\leftarrow} \Y(n).$$
This includes giving a detailed enough description of the action of the automorphism group $\Aut\, \C$ on $\tilde{\HH}$ and, in particular, its action on  $\tilde{\HH}(\mathrm{CM}),$ the special points of the structure (equivalently, the action of $\Gal(\Q(\mathrm{CM}):\Q)$ on $\tilde{\HH}(\mathrm{CM})$), the analogue of (i) of \ref{(i-ii)}. 

It turns out that this problem is closely connected to 
 the 
 \linebreak
theory of   {\em canonical models of Shimura varieties} resulting from Shimura's Conjecture which was developed by G.Shimura, P.Deligne and others. This involves an advanced theory of complex multiplication and Artin's reciprocity map. The results allow one to identify the action of  $\Aut\, \C$ on a single CM-point  (that is the Galois orbit  of the point), see \cite{Milne0}, 12.8. We need a stronger result:  our paper \cite{Append} goes further and identifies the action of $\Gal(\Q(\mathrm{CM}):\Q)$ on $\tilde{\HH}(\mathrm{CM}),$ at the same time describing all the 
relations between CM-points. One of the interesting model-theoretic forms of the main result of \cite{Append} is that {\em the field  $\Q(\mathrm{CM})$ as a structure is $L_{\omega_1,\omega}$-bi-interpretable with a certain structure formulated purely in terms of the ring $\A_f$ of finite adeles over $\Q.$}


In \cite{Append} we defined a certain locally modular structure of trivial type on the set $\tilde{\HH},$ which we term $\tilde{\HH}^\mathrm{Pure}.$ 
This is formulated in terms of the action of a large subgroup $\tilde{\G}$ of
$\GL_2(\A_f)$ on $\tilde{\HH}.$
The main result of \cite{Append} determines automorphisms of $\tilde{\HH}^\mathrm{Pure} $ and then describes  $\Gal(\Q(\mathrm{CM}):\Q)$ in terms of 
$\Aut\,\tilde{\HH}^\mathrm{Pure}(\mathrm{CM}).$ 
These two groups are "almost" equal: there is an obvious embedding of 
$ \Gal(\Q(\mathrm{CM}):\Q)$ into $ \Aut\,\tilde{\HH}^\mathrm{Pure}(\mathrm{CM}),$ and when restricting the two groups to their action on 
a finite number of $\G$-orbits the first is a finite index subgroup of the second.

$\tilde{\HH}^\mathrm{Pure}$  has a ``standard" version $\HH^\mathrm{Pure}$ which is based on the actual upper half-plane $\HH$ and is given in terms of the action of $\GL_2^+(\Q)$ and complex multiplication.

\epk

\bpk With the results of \cite{Append}  in hand we can apply the model-theoretic techniques on categoricity  developed in earlier works (see e.g. the survey in \cite{Zspecial}) and in particular  the proof in \cite{DH} which we follow  quite closely (and so use the Serre open image theorem). 

The axioms $\Sigma$ of the resulting $L_{\omega_1,\omega}$-theory consist, as in abelian examples, of three parts: 

\begin{itemize}
\item[A] The axioms describing $\U^\mathrm{Pure}$ (the abstract version of $\HH^\mathrm{Pure}$),
a locally modular structure of trivial type with the action of the group $\G$ isomorphic to $\GL_2^+(\Q),$ and its subgroups corresponding to $\Gamma(n).$
\item[B] The axioms  describing an algebraically closed field $\F$ of characteristic 0 and
curves $\Y(n)\subset \PP^3(\F).$
\item[C] The axioms describing $j_n: \U^\mathrm{Pure}\to \Y(n)$ obtained via the translation of relations on $\U^\mathrm{Pure}$ to special relations on the $\Y(n).$  
\end{itemize}  

{\bf Main Theorem.} {\em The above system of axioms $\Sigma$ is satisfied by the standard complex structure $\HH$ and every model of the axioms is a quasi-minimal geometry structure.

For each model $\U^\mathrm{Full}$ of $\Sigma$ there
 is a simple $L_{\omega_1,\omega}$-sentence $\Theta$ which holds on $\U^\mathrm{Full}$ and such that $\Sigma\, \&\, \Theta$ has,
  up to isomorphism, a unique model in any uncountable cardinal.}  
\epk
\bpk The first author whishes to express his gratitude to J.Derakhshan for his interest in the work and his help with some of the mathematical issues in the paper.

The second author would like to thank the EPSRC for its support via a New Investigator Award (EP/S029613/1). He would also like to thank the University of Oxford for having him as a Visiting Research Fellow.
\epk

\section{Groups acting on $\HH.$}
\bpk \label{section2} {\bf Groups $\G$ and $\Gamma$ and their generators.} 

Generators of $\Gamma:=\SL_2(\Z)$ are represented by matrices 
$$\sss=\left( \begin{array}{ll}
0 -1\\
1\ \ \ 0
\end{array}\right),\ \  \ttt=\left( \begin{array}{ll}
1\ 1\\
0\ 1
\end{array}\right):$$
$$\sss^2=-\mathbf{I}, \ \ (\sss\ttt)^3=\left( \begin{array}{ll}
0 -1\\
1\ \ \ 1
\end{array}\right)^3=-\mathbf{I}.$$

Generators of $\G:=\GL_2^+(\Q):$

$$\sss,\ttt\mbox{ and } \rrr_q:=\left( \begin{array}{ll}
q\ \ 0 \\
0\ \ 1
\end{array}\right), q\in \Q_+$$
satisfying:
$$\rrr_q\rrr_r=\rrr_{qr},$$
$$ \sss\rrr_q=q\rrr_q\inv \sss,$$
$$\rrr_n \ttt=\ttt^n\rrr_n, \mbox{ for }n\in \N.$$

We will use
$$ \rrr'_q:=\sss  \rrr_{q}\sss\inv=\left( \begin{array}{ll}
1\ \ 0 \\
0\ \ q
\end{array}\right), q\in \Q_+$$ 
$$\ttt_-=\left( \begin{array}{ll}
\ 1\ \  \ 0\\
-1\ \ 1
\end{array}\right),
$$
which satisfies $$\sss \ttt\sss\inv=\ttt_-;\ \ \rrr_n \ttt_-^n=\ttt_-\rrr_n, \mbox{ for }n\in \N,$$
and subgroups
 \be\label{Delta} \Delta(\Q_+)=\{ \rrr_q:\ q\in \Q_+\} \mbox{ and } 
\Delta'(\Q_+)=\{ \rrr'_q:\ q\in \Q_+\}.\ee

Similar notation also make sense for the multiplicative group of a commutative ring $R,$ 
 $$ \Delta(R^\times)=\{ \rrr_q:\ q\in R^\times\} \mbox{ and } 
\Delta'(R^\times)=\{ \rrr'_q:\ q\in R^\times\}.$$

\epk
\bpk \label{AutG} {\bf Remark.} Note that all automorphisms of $\GL_2^+(\Q)$  and of $\mathrm{PGL}_2(\Q)$ are inner and the only non-identity automorphism  which fixes  $\Delta'$ element-wise and preserves the subgroup $\Gamma$ is the involution
$$g\mapsto \check{g}:=\rrr_{-1}\cdot g\cdot \rrr_{-1}.$$

Note also that since the action of the centre of $\GL_2^+(\Q)$ acts trivially on $\HH$ we can work equally with $\mathrm{PGL}_2(\Q),$ the generators and defining relations of which are the same if we 
ignore the scalar multipliers in the group relations. 
\epk

\bpk \label{defspec} {\bf Special points on $\HH.$} Let
$$\E=
 \{ \left( \begin{array}{ll}
a\ b\\
c\ d
\end{array}\right)\in \GG: (d-a)^2<-4bc\}= \{ e\in \GG: \mathrm{tr}^2\,e<4\det e\},$$
the set of {\em elliptic} transformations.

These are exactly the elements for which
there is a unique fixed point $\tau_e\in \HH$ which
can be found by solving the equation \be \label{ch} c x^2+(d-a)x-b=0,\ \ \tau_e:=x,\ \ \Im x>0.\ee

Note that elements of the centre 
$$\mathrm{Z}=\{ \left( \begin{array}{ll}
a\ 0\\
0\ a
\end{array}\right)\in \G\}$$
act as identity on $\HH.$
 \epk
\bpk \label{defE}
{\bf Remark.} (i) $\E$ is 
 invariant under automorphisms of $\G$ since  $g\cdot e\cdot g\inv$ fixes $g\tau_e$ and thus belongs to $\E.$

(ii) The subgroup $\mathrm{St}_\tau\subset\G$ fixing point $\tau=\tau_e$ is definable by the condition  that   $\mathrm{St}_\tau=C(e),$ the centraliser of $e$ in the group. This is a general fact for Mobius transformations over any field of characteristic 0.

\epk
The following is a technical remark which we are going to use later.
\bpk \label{listG} {\bf Remark.} Let $I=\bigcup_l \E^l$ and $\{ \Phi_h: h\in I\}$ a family  of quantifier-free $L_{\omega_1,\omega}$-formulas.

 There is an existential $L_{\omega_1,\omega}$-formula $\Theta$ 
 in the   language with names $\rrr_q$ for respective elements of $\G$ 
  stating that
for any $\mathsf{h}\in \E^l$ there exists $\bar{t}_h=\la t_1,\ldots,t_l\ra,$ a tuple of fixed points of $\mathsf{h}$ such that $\Phi_i(\bar{t}_h).$ 

Indeed, let $\{ w_e(\sss,\ttt, \rrr_\Q):\ e\in \E\}$ be the family of group words in generators $\sss, \ttt,\rrr_q$ $q\in \Q_+,$  listing all the elements of $\E.$ Using the family of words we can  produce the family  $\{ \bar{w}_h(\sss,\ttt,\rrr_\Q):\ \mathsf{h}\in I\}$ of tuples of words corresponding to tuples $\mathsf{h}\in I.$ Let $\mathrm{Fix}(\bar{t}, \bar{h})$ is the formula saying that  $\bar{t}$ is a tuple from $\HH$ fixed by the tuple $\bar{h}$ from $\E.$

Set
$$\Theta:= \begin{array}{ll}
 \exists s,t\in \G\ \{ s^2=-1\ \& & (st)^3=-1\ \&  \ \bigwedge_{q\in \Q_+} \sss\rrr_q=q\rrr_q\inv \sss \ \&\ \bigwedge_{n\in \N} \rrr_n\ttt=\ttt^n\rrr_n \ \& \\ 
&\&\ \bigwedge_{  \mathsf{h}\in I} \exists \bar{t}_h \mathrm{Fix}(\bar{t}_h, \bar{w}_h(s,t,\rrr_\Q))\ \& \ \Phi_h(\bar{t}_h)\}\end{array}$$
 
Since the choice of any $s,t\in \G$ satisfying the group relations in the first line of $\Theta$ is conjugated to $\sss,\ttt$ by an automorphism of $\G$ (see  \ref{AutG}) and $\E$ is invariant under the automorphisms, the values of words $\bar{w}_h(s,t,\rrr_\Q))$ run through all elements of $\E.$
 $\Box$  
\epk
\section{Projective limit}\label{section3}
\bpk \label{pr-limit}  We consider the structure $\tilde{\HH},$ the projective limit of structures
\be\label{Htilde} \tilde{\HH}=\lim_{\leftarrow }\Gamma
(N)\backslash \HH,\ee
along the projective system of 
$$\Gamma(N)=\left\lbrace \left( \begin{array}{ll}
a\ b\\
c\ d
\end{array}\right)\in \Gamma,\ \ a\equiv d\equiv 1\mod N,\ 
c\equiv b\equiv 0\mod N\right\rbrace,
$$ normal subgroups of $\Gamma=\mathrm{SL}(2,\Z).$ 

The structures  $\Gamma(N)\backslash \HH$ are identified with the classical complex modular curves $\Y(N)$ (without the full level $N$ structure). These are non-singular algebraic curves and therefore  can be realised as  quasi-projective algebraic curves in $\PP^3(\C).$ Moreover, the realisation can be obtained over  $\Q,$ see e.g. section 4 of \cite{Elkies} or the paper \cite{Xavier}, section 2. Another argument for this fact is given in \cite{Append}, 3.2(8).
 
 The classical modular functions 
$$j_N:\HH\to \Y(N);\ \ \tau\mapsto \Gamma(N)\cdot\tau$$ 
are holomorphic
and the curves $\Y(N)$ are finite
coverings of  $\Y(N/d),$ for $d|N,$  via the projection maps
\be\label{covers} \pr_{N, N/d}: j_N(\tau)\mapsto j_{N/d}(\tau),\ \ \ \Y(N)\to \Y(N/d).\ee

\epk
\bpk The definition (\ref{Htilde}) implies the existence of maps
$${j}_N: \tilde{\HH}\to \Gamma(N)\backslash \HH=\Y(N)$$
which we are going to use below.

By definition, any $\tau\in \tilde{\HH}$ is uniquely determined by 
the sequence $${j}_n(\tau)\in \Y(n), \ n\in \N$$
and the sequence has the property
$$\pr_{n, n/d}({j}_n(\tau))={j}_{n/d}(\tau), \mbox{ for each } d|n.$$
\epk 
 \bpk\label{3.7} {\bf Lemma.} {\em  Any $\tau\in \HH$ is uniquely determined by the sequence $${j}_n(\tau)\in \Y(n): \ n\in \N.$$ This gives  the canonical embedding $$\HH\hookrightarrow \tilde{\HH}.$$  }

Proof. The first statement follows from the fact that $\bigcap_N \Gamma(N)=\{1\}.$
The second statement is the consequence of the fact that the sequence satisfies
$$\pr_{n, n/d}(j_n(\tau))=j_{n/d}(\tau), \mbox{ for each } d|n.$$
 $\Box$
\epk

\bpk \label{general} {\bf Remark.} The system of covers (\ref{covers}) is not \'etale. However, by removing finitely many  points on $\Y(1)$ and all the points on $\Y(n)$ over these one gets smooth curves $\Y^-(n)$ and a projective system of covers $$\pr_{n,n/d}: \Y^-(n)\to \Y^-(n/d)$$
  which  is \'etale. 

Since the construction of $\tilde{\Gamma}$ depends on generic fibres we have the same  $\tilde{\Gamma}$ for the construction corresponding to the system of \'etale covers. 
\epk
\bpk \label{Hfull} {\bf The projective limit as the structure $\tilde{\HH}^\mathrm{Full}.$}

The analysis and study of the projective limit of modular curves $\Y(n)$ with the {\em $n$-level structure defined over $\Q(\zeta_n)$} is readily reducible to the study of {\em canonical models of Shimura varieties}, see \cite{Milne0}, and more specifically   canonical models of modular curves,  \cite{Milne1}. The case of curves $\Y(n)$  over $\Q$ as above in the setting appropriate for our purposes is studied in \cite{Append}.
The structure $\tilde{\HH}^\mathrm{Full}$ is described therein as the quotient $\Delta(\hat{\Z}^\times)\backslash\mathbf{S}^\mathrm{Full}$ of a more fundamental structure $\mathbf{S}^\mathrm{Full}$ associated with the Shimura datum $(\HH\cup -\HH, \GL_2).$

\medskip

The main conclusions, see \cite{Append}, 3.22, are as follows:

(a) There is a group $\tilde{\G}$ acting on $\tilde{\HH}.$ 
$$\tilde{\G}\cong \Delta(\Q_+)\cdot\SL_2(\A_f)\subset \GL_2(\A_f),$$
where $\A_f$ is the ring of finite adeles.

One can speak about an action of $\gh\in \Delta(\Q_+)\cdot\SL_2(\A_f)$ on $\tilde{\HH}$ once an isomorphism $\varphi: \Delta(\Q_+)\cdot\SL_2(\A_f)\to \tilde{\G}$ is provided. Call $\varphi$ a {\bf naming isomorphism}.   The naming isomorphisms  form the family
$$\{ \varphi_\mu:\  \Delta(\Q_+)\cdot\SL_2(\A_f)\to \tilde{\G} \ |\ \mu \in 
\hat{\Z}^\times\};\ \mbox{ set } \gh^\mu:= \varphi_\mu(\gh),$$
which satisfies 
$$ \gh^{\mu\cdot \lambda}=(\rrr_\lambda\inv\cdot \gh\cdot \rrr_\lambda)^\mu$$
where the conjugation by $\rrr_\lambda$ is in the ambient group $\GL_2(\A_f).$

\medskip

(b) 
Any $\gh\in \Delta(\Q_+)\cdot\SL_2(\A_f)$ determines 
the 0-definable subset of $\tilde{\G}:$
$$\gh^\Delta=\{ \gh^\mu: \mu\in \hat{\Z}^\times\}.$$

Each $\gh^\mu$ gives rise to
 the sequence of algebraic curves defined over $\bar{\Q},$
$$C^\mu_{\gh,N}\subset \Y(N)\times \Y(N);\ \ C^\mu_{\gh,N}:=\{ \la {j}_N(u),{j}_N(\gh^\mu\cdot u)\ra: u\in \tilde{\HH}\} $$
(here $C^\mu_{\gh,N}$ corresponds to $C^\mu_{\gh,\approx K}$ of \cite{Append}). Note that any special curve $C\subset \Y(N)\times \Y(N)$ has the form $C=C^\mu_{\gh,N}$ for some $\gh.$

The construction of the projective limit applied to the sequence $\{ C^\mu_{\gh,N}: N\in \N\}$ of curves results in the {\em limit curve }
$$C^\mu_{\gh}\subset \tilde{\HH}\times \tilde{\HH};\ \ C^\mu_{\gh}:=\{ \la u,\gh^\mu u\ra: u\in\tilde{\HH}\}.$$

(c) For a fixed $\gh\in \Delta(\Q_+)\cdot\SL_2(\A_f)$ we obtain the finite family of curves on $\Y(N)\times \Y(N)$   $$\{ C^\mu_{ \gh,N}, \ \mu\in \hat{\Z}^\times\}$$
where $$\mu-\lambda\in N\hat{\Z}\Longrightarrow  C^\mu_{\gh,N}= C^\lambda_{\gh,N} .$$ 
These curves are irreducible components of the algebraic curve $C_{\gh,N}$ defined over $\Q:$ 
$$C_{\gh,N}= \bigcup_{\mu\in \hat{\Z}^\times} C^\mu_{\gh,N}\subset \Y(N)\times \Y(N).$$

In the limit one obtains the the  infinite-component limit curve on $\tilde{\HH}\times \tilde{\HH}:$
$$C_\gh= \bigcup_{\mu\in \hat{\Z}^\times} C^\mu_{\gh}.$$
$C_\gh$ is defined over $\Q$ too.

(d) The irreducible components $C^\mu_{\gh,N}$ of $C_{\gh,N}$ are Galois conjugated over $\Q.$

Clearly $C_\gh=C_{\gh'}$ for $\gh'=\gh^\mu.$
The image of $C_\gh$ under $j_N\times j_N$ is $C_{\gh,N}.$

\medskip

(e) For any $\gh\in  \Delta(\Q_+)\cdot\SL_2(\A_f),$
the definable 4-ary relation on $\tilde{\HH}$  
$$\mathrm{Comp}_\gh(s_1,s_2,t_1,t_2):\equiv \exists h\in \tilde{\G}\
s_2=h\cdot s_1\ \& \ t_2=h\cdot t_1\ \&\ C_\gh(s_1,s_2)$$
determines the condition that $\la s_1,s_2\ra$ and $\la t_1,t_2\ra$ belong to the same irreducible component of a $C_\gh.$  

For each $N$ the relation $\mathrm{Comp}_{\gh,N}$ on $\Y(N)$ is the image of $\mathrm{Comp}_\gh$ under $j_N,$ the relation which determines the decomposition of the algebraic curve $C_{\gh,N}$ into irreducible components.  This relation is invariant under $\Gal_\Q$ and so definable over $\Q.$

 We also have on $\tilde{\HH}$ the definable equivalence 
$$j_N(s_1)=j_N(s_2),$$
which can be equivalently given by:
$$\exists \gamma\in \tilde{\Gamma}(n)\ s_2=\gamma\cdot s_1.$$

\medskip 

(f) The points $s$ in $\tilde{\HH}$ which are fixed by a $g\in \tilde{\G}\setminus \varphi_\mu(\mathrm{Z})$ (the centre of $\GL_2^+(\Q)$
 will be called special, or CM-points of $\tilde{\HH}.$ $j_N(s)$ is a CM-point in $\Y(N)$ and, in particular, is algebraic.

The paper \cite{Append} describes the binary relations on $\tilde{\HH}$ written as
$$\la t_1,t_2\ra\in \tp(s_1,s_2)$$
defined for each pair  $s_1,s_2$ of CM-points. The reation is valid if and only if there is an automorphism $\sigma$ of the projective system (\ref{covers}) such that $\sigma: \la s_1,s_2\ra \mapsto \la t_1,t_2\ra.$ The relations are 
invariant under automorphisms  of the projective system (i.e. defined over $\Q$) for each choice of $s_1,s_2.$

 The image of the relation under $j_N$ is the relation on $\Y(N)$ 
$$\la y_1,y_2\ra\in\tp_N(x_1,x_2)$$
which holds if and only if $\la y_1,y_2\ra$ is Galois conjugated to  $\la x_1,x_2\ra$ over $\Q.$


\epk
\bpk {\bf Definition.} $\tilde{\HH}^\mathrm{Pure}$ is the structure with the universe $\tilde{\HH}$ and relations $C_\gh(s_1,s_2),$ 
$\mathrm{Comp}_\gh(s_1,s_2,t_1,t_2)$ ($\gh\in \Delta(\Q_+)\cdot\SL_2(\A_f)$), 
$j_N(s_1)=j_N(s_2)$ and $\la t_1,t_2\ra\in\tp(s_1,s_2).$

\medskip

$\tilde{\HH}^\mathrm{Full}$ is the multisorted structure with sorts
$\tilde{\HH}$ and $\Y(N),$ $N\in \N,$ and relations: 

- on   $\tilde{\HH}$ the relations of  $\tilde{\HH}^\mathrm{Pure};$

- on the $\Y(N)$ the Zariski closed relations defined over $\Q;$

- the maps $j_N:  \tilde{\HH}\to \Y(N)$ and $\pr_{N, N/d}:  \Y(N)\to \Y(N/d).$

\medskip

$\tilde{\HH}^\mathrm{Pure}(\mathrm{CM})$ and $\tilde{\HH}^\mathrm{Full}(\mathrm{CM})$ are substructures of the structures with universes restricted to their special points. 

\medskip

{\bf Remarks.} (i) Since components of the $C_\gh(s_1,s_2)$ are graphs of actions of elements of the group $\tilde{\G},$ we often look at $\tilde{\HH}$ as a  $\tilde{\G}$-set. 

(ii) A corollary of definitions is that the relations $C_\gh,$ 
$\mathrm{Comp}_\gh,$  and $\la t_1,t_2\ra\in\tp(s_1,s_2)$
are projective limits of the relations $C_{\gh,N},$ 
$\mathrm{Comp}_{\gh,N},$  and $\la t_1,t_2\ra\in\tp_N(s_1,s_2)$ on the $\Y(N).$ That is the relations on the sort $\tilde{\HH}$ are positive-type-definable in terms of pull-backs of the respective relations on the $\Y(N)$ along with the pull-back of equality.  

We thus can, up to $L_{\omega_1,\omega}$-bi-interpretability, identify  $\tilde{\HH}^\mathrm{Pure}$ as the structure given by the pull-backs of $C_{\gh,N},$ 
$\mathrm{Comp}_{\gh,N},$   $\la t_1,t_2\ra\in\tp_N(s_1,s_2)$ and equality.
Call it {\bf the pull-backs version of  $\tilde{\HH}^\mathrm{Pure}.$} In the same sense we speak about the pull-backs version of  ${\HH}^\mathrm{Pure},$ a substructure of the pull-backs version of  $\tilde{\HH}^\mathrm{Pure}.$

\epk
\bpk {\bf Lemma.} {\em The  pull-backs versions of  $\tilde{\HH}^\mathrm{Pure}$ and of ${\HH}^\mathrm{Pure}$ are locally modular of trivial type.

The structures satisfy
$${\HH}^\mathrm{Pure} \prec \tilde{\HH}^\mathrm{Pure}.$$  }

{\bf Proof.} The structure ${\HH}^\mathrm{Pure}$ satisfies the assumptions of Theorems 4.11 and 4.14 of \cite{Zspecial}. It follows that its theory has quantifier elimination and is of trivial type. $\tilde{\HH}^\mathrm{Pure}$ is then its elementary extension in the obvious way.  
 $\Box$

\epk
\bpk {\bf Remarks.} It is easy to see that the centre $\mathrm{Z}$ of $\GL^+_2(\Q)$ 
acts on $\tilde{\HH}$ trivially.

\medskip

Note that in the definitions above the curves $\Y(N)$ are curves over $\C$ 
(defined over $\Q$) and so the points of $\tilde{\HH}$ are limits of $\C$-points. 

However the definitions and results are valid in the context of curves $\Y(N)$   over $\F,$ an abtract algebraically closed field of characteristic $0.$ In this case $\Y_\F(N)$ are curves over $\F.$ 
\epk
\bpk 
Define $\tilde{\U}_\F^\mathrm{Pure}$ and $\tilde{\U}_\F^\mathrm{Full}$ to be the respective structures obtained as the projective limit of
the $\Y_\F(N),$ for an arbitrary algebraically closed field of characteristic $0.$
  
\epk

\bpk \label{U(CM)}
{\bf Remark.} Note that $$
\tilde{\HH}^\mathrm{Pure}(\mathrm{CM})=\tilde{\U}_\F^\mathrm{Pure}(\mathrm{CM})$$
since both sides are the structures obtained by taking the projective limit of the respective substructures  $\Y(n)^\mathrm{Pure}(\mathrm{CM})$ on the curves. 
The same is true for the full  structures:
$$\tilde{\HH}^\mathrm{Full}(\mathrm{CM})=\tilde{\U}_\F^\mathrm{Full}(\mathrm{CM}).$$
\epk

\bpk\label{deftE} Define
$$\tilde{\E}=\{ e\in \tilde{\G}\setminus \mathrm{Z} \ \exists u\in \tilde{\U}\ e\cdot u=u\}.$$


Set

$$\tilde{\G}_*=( \tilde{\G}, \tilde{\Gamma}, \tilde{\E}, \{ \rrr_q, \ \rrr'_q: q\in \Q_+)$$
the structure on the group $\tilde{\G}$ with distinguished subgroup
$\tilde{\Gamma},$ distinguished subset $\tilde{\E}$ and distinguished elements $\rrr_q$ and  $\rrr'_q.$

Analogously,
$${\G}_*=( {\G}, {\Gamma}, {\E}, \{ \rrr_q, \ \rrr'_q: q\in \Q_+\}).$$

Clearly, $$\G_*\subset  \tilde{\G}_*$$
as structures.  
\epk
\bpk {\bf Remark.} There is an embeding 
$$\tilde{\G}\subset \GL_2(\A_f)$$
which is an identity on the diagonal elements  $\rrr_q, \ \rrr'_q: q\in \Q_+.$ It is easy to see that such an embedding is determined uniquely, up to the conjugation by $\rrr_\mu,$ $\mu\in \hat{\Z}^\times.$

In particular, we may assume that elements $g$ of $\tilde{\G}$ are also elements of $\GL_2(\A_f)$ and thus the conjugation by an element $\rrr_\lambda\in \Delta,$
$$g\mapsto \rrr_\lambda\cdot g\cdot \rrr_\lambda\inv$$
is well defined.
\epk

\bpk\label{Claim1} {\bf Lemma.} {\em Consider the natural embedding $\G_*\subset \tilde{\G}_*$ and let 
$$\psi: \G_* \twoheadrightarrow \G'_*\subset \tilde{\G}_*$$
be a partial isomorphism of  $\tilde{\G}_*.$ 

$\psi$ can be extended to an automorphism $\tilde{\psi}: \tilde{\G}_*\to \tilde{\G}_*.$ Moreover, there is $\lambda\in \hat{\Z}^\times$ such that
$$\tilde{\psi}:g\mapsto \rrr_\lambda\cdot g\cdot \rrr_\lambda\inv,\mbox{ for all } g\in \tilde{\G}.$$
}

{\bf Proof.} Let $\sss'=\psi(\sss)\in \G'.$ Then
$$\sss'\in \tilde{\G}, \ \ \sss'\cdot \rrr_{-1}=\rrr'_{-1}\sss'\mbox{ and }
\ \ {\sss'}^2=-\mathbf{I}=\rrr_{-1}\cdot\rrr'_{-1}.$$
It is easy to see that the three equations imply that, for some $\lambda\in \hat{\Z}^\times,$
$$\sss'=\left(\begin{array}{ll}
\ 0\ \ \ \ \lambda\\ -\lambda\inv\ 0
\end{array}\right)=   \rrr_\lambda\cdot \sss\cdot \rrr_\lambda\inv.$$

Let $\ttt'=\psi(\ttt)\in \G'.$ Then 
  for each  $n\in \N_{>0},$ 
$$ \ttt'\in \tilde{\Gamma},\ \  \rrr_n\cdot \ttt'\cdot \rrr_n\inv={\ttt'}^n \mbox{ and } (\sss' \ttt')^3=\mathbf{I}.$$
It follows 
$$\ttt'=\left(\begin{array}{ll}
1\  \lambda\\ 0\ 1
\end{array}\right)=   \rrr_\lambda\cdot \ttt\cdot \rrr_\lambda\inv.$$
Thus, $\G'$ is the group generated by $\sss'$ and $\ttt'$ and 
$$\psi: g\mapsto \rrr_\lambda\cdot g\cdot \rrr_\lambda\inv\mbox{ for all } g\in \G.$$
Take $$\tilde{\psi}:g\mapsto \rrr_\lambda\cdot g\cdot \rrr_\lambda\inv\mbox{ for all } g\in \tilde{\G}.$$

It is clear that $\tilde{\psi}$ preserves $\tilde{\Gamma}$, $\rrr_q$ and $\rrr'_q,$ $q\in \Q_+.$ 

Finally note that 
$g\in \tilde{\E}$ if and only if 
$\tilde{\psi}(g)\in \tilde{\E}.$ This follows from the description of fixed points in  3.22D(b) of \cite{Append}. 

$\Box$
\epk


The following is essentially a corollary of main results of \cite{Append}.

\bpk \label{Claim2} {\bf Lemma.} {\em
 Any  $\psi\in \Aut\,  \tilde{\G}_*$ can be extended to $\tilde{\psi}\in \Aut\, \tilde{\U}_\F^\mathrm{Pure}(\mathrm{CM}).$
 }
 
 {\bf Proof}  As noted above  $\tilde{\U}_\F^\mathrm{Pure}(\mathrm{CM})=\tilde{\HH}^\mathrm{Pure}(\mathrm{CM})$ so we may argue in the setting of $\tilde{\HH}.$

  By \ref{Claim1} $\psi$ has the form 
  $\gh^\mu\mapsto \gh^{\mu\cdot \lambda}= \rrr_\lambda\cdot \gh^\mu\cdot \rrr_\lambda\inv,$ for some $\lambda\in \hat{\Z}^\times.$ And by \ref{Hfull}(d) $\psi$ is induced by a $\sigma\in \Gal_\Q.$ In its turn $\sigma$ acts on the CM-points of $\tilde{\U}_F$ and induces a $\tilde{\psi}\in \Aut\, \tilde{\U}_\F^\mathrm{Pure}(\mathrm{CM}),$ as required.
  $\Box$

\epk
\section{Axiomatisation of $\HH_\F^\mathrm{Full}$ and $\U_\F^\mathrm{Full}.$} \label{section 5}

\bpk {\bf Axioms}.

The language $\LL(j_{ n\in \N})$ is 3-sorted, with sorts $\U,\G$ and $\F.$
The structure on $\F$ is that of a field given in a standard ring language, the structure on $\G$ and $\U$ is that of a group acting on $\U$ with distinguished subsets $\E, \Gamma, \{\rrr_q,  \rrr'_q, q\in \Q^\times\}.$ Note that  $\{\rrr_q, q\in \Q^\times\}=\Delta(\Q)$ with all its elements named and the same for
$\{\rrr'_q, q\in \Q^\times\}=\Delta'(\Q).$
Note that $$\GL_2^+(\Q)=\Delta(\Q_{>0}) \cdot \SL_2(\Z)=\Delta'(\Q_{>0}) \cdot \SL_2(\Z)$$ and so  group $\G$ is isomorphic to $\GL_2^+(\Q)$ with distinguished elements of $\Delta,$ $\Delta'$ and subgroup $\Gamma\cong \SL_2(\Z)$ will have the same structure. 

Note that   $\GL^+_2(\Q)$ is invariant under the involutive transformation
 $$g\mapsto \check{g}:=\rrr_{-1}\cdot g\cdot \rrr_{-1}$$
 where $\rrr_{-1}\in \GL_2(\Q)\subset \GL_2(\A_f).$
 
    The maps $j_n$ have $\U$ as their domain and have quasi-projective curves $\Y(n)$ as their range.  
   $\Sigma$ consists of the following five groups of axioms :

\medskip

Group  axioms: \be \label{GG} (\G,\Gamma,\E, \{ \rrr_q,\rrr'_q: q\in \Q_{>0}\})\cong (\GL^+_2(\Q), \SL_2(\Z), \E(\Q), \{ \rrr_q,\rrr'_q: q\in \Q_{>0}\}).\ee

Note that $\Gamma(n)$ is definable from the data, namely the subgroup $\Gamma_0(n)$ of
matrices of the form $\gh=\ttt^{nm}$ are definable by condition
$$\exists \gamma\in \Gamma\ \gh=\rrr_n\cdot \gamma\cdot \rrr_n\inv,$$
and $\Gamma(n)$ can be defined as the normal closure of  $\Gamma_0(n).$

 Let $\mathrm{Z}=\mathrm{Centre}(\G).$

\medskip

{\bf Action axiom}:  $\G \mbox{ acts on } \U;$
\be \label{AA} \begin{array}{lll}
 \forall g\in \G\setminus (\E\cup \mathrm{Z})\ \forall u\in \U\
 g\cdot u\neq u,\\

\forall g\in \mathrm{Z}\
\forall u\in \U\ g\cdot u=u,\\

\forall e\in  \E\
\exists! u_e\in \U\ e\cdot u_e=u_e. 
 \end{array}
  \ee
  
 {\bf Fibre formula}:

\be \label{Ff}
\forall u,v\in \U\ j_n(u)=j_n(v)\leftrightarrow \exists \gamma\in \Gamma(n)\ v=\gamma\cdot u\ee 

\medskip

{\bf ACF$_0$ axioms and sorts $\Y(n)$}:

\be \label{ACF} \F\vDash \mathrm{ACF}_0
\ee 
and $$\Y(n)\subset \PP^3(\F); \ \ \pr_{n.m}: \Y(n)\to \Y(m),\mbox{ for }m|n$$  are given by specific equations over $\Q.$
\medskip

{\bf Functional equations}: 
\be \label{A01} j_n: \U\twoheadrightarrow \Y(n); \ \ \pr_{n,m}\circ j_n=j_m\mbox{ for each }m|n;
\ee

\be \label{A2}\begin{array}{ll} \forall g\in \G\ \exists \gh\in \GL^+_2(\Q), \exists \mu\in \hat{\Z}^\times: \ \bigwedge_{n\in \N}   C^\mu_{\gh,n}=j_n(\mathrm{graph}\,g)\\
\forall q\in \Q_+\ \forall \mu\in \hat{\Z}^\times: \ \bigwedge_{n\in \N}   C^\mu_{\rrr_q,n}=j_n(\mathrm{graph}\,\rrr_q)
\end{array}\ee

 $\forall \gh\in \GL_2^+(\Q),$  $\forall \mu\in \hat{\Z}^\times,$  $\forall u,v\in \U:$  
\be \label{B1}   \bigwedge_{n\in \N} \la j_n(u),j_n(v)\ra\in C^\mu_{\gh,n} \Leftrightarrow \exists g\in \G\  v=g\cdot u\ \&\ \bigwedge_{n\in \N} C^\mu_{\gh,n}=j_n(\mathrm{graph}\,g)\ee

\medskip
\epk

\bpk {\bf Proposition.} {\em $\HH_\C^\mathrm{Full}$ is a model of $\Sigma.$}

{\bf Proof.} Axioms (\ref{GG})-(\ref{A2}) just  list general properties of $\HH$ and $j_n$ from section~\ref{section3} (axioms (\ref{A01})-(\ref{A2})   are established in \ref{Hfull}).

Axiom (\ref{B1}) is proved for $\HH$ in \cite{Append} in the Claim 4.8(5). Indeed, the left-hand side of (\ref{B1}) states that $\la u,v\ra$ is a point on a graph of an element $\gh^\mu\in \tilde{\G}$ whose graph in $\tilde{\HH}$ is $C^\mu_\gh.$ By the Claim $\gh^\mu\in \GL^+_2(\Q).$

 $\Box$
\epk
\bpk\label{subs} {\bf Lemma.} {\em For any model $\U_\F^\mathrm{Full}$ of $\Sigma$ 
 there is an embedding
\be\label{subsi} \mathbf{i}: \U_\F^\mathrm{Full}\hookrightarrow \tilde{\U}_\F^\mathrm{Full}\ee
which is the identity on  $\F$ and induces an embedding $\mathbf{i}:\G/Z\hookrightarrow \tilde{\G}/\{ \pm 1\}$ such that the $\rrr_q$-named elements  of $\G$ map to the $\rrr_q$-named elements of $\tilde{\G},$ $\mathbf{i}(\Gamma)\subset \tilde{\Gamma}$ and  $\mathbf{i}(\E)\subset \tilde{\E}.$

}

{\bf Proof.} $\tilde{\U}$ in $\tilde{\U}_\F^\mathrm{Full}$ is given as the universal cover of $\Y(\F),$ with respect to the class of modular curves. By axiom (\ref{A01})    $\U$  covers every principal modular curve $\Y(n)$ by $j_n: \U\to \Y(n)$ which agrees with the system $\pr_{n,m}$ by the same rule as $\tilde{j}_n$ does. Hence, for any  $u\in \U$ there is a unique $\tilde{u}\in \tilde{\U}$ such that
$j_n(u)=\tilde{j}_n(\tilde{u})$ for all $n.$
 This determines a unique  map $$\mathbf{i}: u\mapsto \tilde{u};\ \U \to \tilde{\U}.$$

Axioms (\ref{GG}) and (\ref{Ff}) ensure that $\mathbf{i}$ is an embedding. (\ref{ACF}) and (\ref{A01}) imply that the structure on the $\Y(n)$ in $\U_\F$ and in $\tilde{\U}_\F$ are identical and that the maps $j_n$ and $\tilde{j}_n$ satisfy the same equations.
The first line of (\ref{A2}) ensures that $\mathrm{i}$  sends the graph of $g\in \G$ to the relations $C^\mu_{\gh,n}$ for all $n,$ equivalently $\mathrm{graph}\, g$ corresponds to the graph of a $\gh^\mu\in \tilde{\G}.$  Thus $\mathbf{i}: g/Z\mapsto \gh^\mu/Z$  is  an embedding  of $\G/Z$  into $\tilde{\G}/(Z\cap \tilde{\G})=\tilde{\G}/\{ \pm 1\}$ such that
$\mathbf{i}(g*u)=\gh^\mu * \mathbf{i}(u)= \mathbf{i}(g)*\mathbf{i}(u).$
The second line of (\ref{A2}) tells us that $\mathbf{i}(\rrr_q)=\rrr_q$
(note that $\rrr_q^\mu=\rrr_q,$ for all $\mu$).

We get  $\mathbf{i}(\Gamma(n))\subset \tilde{\Gamma}(n)$ by (\ref{Ff}),
and $\mathbf{i}(\E)\subset \tilde{\E}$ by (\ref{GG}).

$\Box$

\epk

\bpk\label{gg'} {\bf Lemma.} {\em Suppose  $\gh\in \GL_2^+(\Q),$ $g,g'\in \G$ and for some $\mu, \mu'\in \hat{\Z}^\times,$  in structure $\tilde{\U}:$
\be\label{eqgg'} \bigwedge_{n\in \N} C^\mu_{\gh,n}=j_n(\mathrm{graph}\,g) \mbox{ and } \bigwedge_{n\in \N} C^{\mu'}_{\gh,n}=j_n(\mathrm{graph}\,g')\ee

Then $g'=g$ or $g'=\hat{g},$ where
$\hat{g}:=\rrr_{-1} g \rrr_{-1}.$ 

Moreover, there is $\sigma\in \Gal_\Q$ such that $g'=g^\sigma.$ 
   }
  
{\bf Proof.} The equalities imply that $C^\mu_\gh=\mathrm{graph}\,g$ and
$C^{\mu'}_\gh=\mathrm{graph}\,g',$ that is in the imbedding $\G\subset \tilde{\G}$ 
$$g=\gh^\mu\mbox{ and } g'=\gh^{\mu'}$$
that is $g'=g^\lambda,$ for  $\lambda\in \hat{\Z}^\times$ defining an
  automorphism $\psi_\lambda: \tilde{\G}\to \tilde{\G},$ $g\mapsto g^\lambda.$ 

We may assume $\G=\GL_2^+(\Q)\subset \tilde{\G}.$   
  If $g\neq g'$ then $g=\left( \begin{array}{ll}
  a\ b\\ c\ d
  \end{array}\right)$ is not a diagonal matrix, say $b\neq 0,$ and  $g'=g^\lambda=\left( \begin{array}{ll}
  a\ \ \ \lambda b\\ \lambda\inv c\ d
  \end{array}\right).$ Since $b$ and $\lambda b$ are distinct non-zero rational numbers, we have necessarily $\lambda=-1.$

Finally note that $\psi_\lambda$ is an automorphism of $\tilde{\G}$ which is induced by a Galois automorphism acting on $\tilde{\U}_\F^\mathrm{Full}.$    
  $\Box$

\epk

\section{An analytic Zariski structure}
\bpk Let $\kk\subseteq \F$ be a subfield and $S\subseteq \U^m$ in a model $\U_\F$ of $\Sigma.$ We say $S$ is a {\bf closed} (analytic) set over $\kk$ if there is a family $\{ W_n\subseteq \Y(n)^m: n\in \N\}$ of Zariski closed subsets defined over $\kk$ such that
$$S=\bigcap_{n\in \N} j_n\inv(W_n).$$

$S$ over $\kk$ is said to be {\bf irreducible} over $\kk'$ ($\kk\subseteq \kk'\subseteq \F$) if for any countable family
$\{ S_i\subseteq \U^m: i\in I\}$ of  closed sets over $\kk'$ 
$$S=\bigcup_{i\in I} S_i \Rightarrow S=S_{i_0},\mbox{ for some }i_0\in I.$$

For $\bar{u}\in \U^m$ we call $\mathrm{locus}(\bar{u}/\kk)$ the smallest
closed $S\subseteq \U^m$ over $\kk$ which contains $\bar{u}.$ This is, clearly, irreducible over $\kk.$
 
We say that $\bar{u}=\la u_1,\ldots,u_m\ra\in \U^m$ is $\G$-generic if there is no $g\in \G$ such that $u_j=g\cdot u_i$ for $1\le i,j\le m.$ 
\epk
In \ref{L6.2}-\ref{P6.7} we prove, as the matter of fact, that the universe $\U$ equipped with predicates for analytic sets over $\Q$ is a one-dimensional {\em analytic Zariski structure} as defined in \cite{Zbook}. Such a structure, as proved  in   \cite{ZAnZar} (see also \cite{Misha}), is a quasi-minimal geometry structure which, according to \cite{K5},  can be axiomatised categorically in uncountable cardinals by an $L_{\omega_1,\omega}(Q)$-sentence. The is summarised in the 
 final Theorem \ref{MThm}. 
\bpk\label{L6.2} {\bf Lemma.} {\em Let $\bar{u}\in \U^k$ be $\G$-generic. Then there is an $n\in \N$
and a Zariski closed $Z_{\bar{u},\Q}\subseteq \Y(n)^k$ over $\Q$ satisfying
$$\mathrm{locus}(\bar{u}/\Q)=j_n\inv(Z_{\bar{u}},\Q).$$

Let $\F_0\subset \F$ be an algebraically closed subfield and assume that at least one coordinate of $j(\bar{u})$ is not in $\F_0.$ Then 
 there is an $n\in \N$
and a Zariski closed $Z_{\bar{u},\F_0}\subseteq \Y(n)^k$ over $\F_0$ satisfying
$$\mathrm{locus}(\bar{u}/\F_0)=j_n\inv(Z_{\bar{u},\F_0}).$$

}

{\bf Proof.} We use \cite{DH}, section 5.1.  
 Let $\kk_0=\Q^\mathrm{ab}(\mathrm{CM}(1)),$  the extension of $\Q^\mathrm{ab}$ by the co-ordinates of special points in $\Y(1).$ Let $\bar{a}=j(\bar{u}).$ 

$\Gal(\F/\kk_0(\bar{a}))$ acts on $\tilde{\U}$ so that 
 there is a number $M\in \N$ and a  subgroup $\Omega\subseteq \tilde{\Gamma}^k$ of index $M,$ $$\tilde{\Gamma}^k= 
 \bigcup_{1\le i\le M} \Omega\cdot\bar{\gamma}_i, \ \ \ 
  j\inv(\bar{a})=\dot{\bigcup_{1\le i\le M}} \Omega \cdot\bar{u}_i, \ \mbox{ where } \bar{u}_i:=\bar{\gamma}_i \cdot\bar{u}$$ \ \  $$\Omega\cdot\bar{u}_i=\Gal(\F/\kk_0)\cdot\bar{u}_i, \mbox{ for }i=1,\ldots, M.$$ 
 Assume $\bar{u}=\bar{u}_1.$ 
 
Each of the $M$ Galois orbits correspond to the algebraic type of the
sequence $\{ j_n(\bar{u}_i): n\in \N\}$ given by systems of equations
\be \label{eqP}\bar{P}_{i,n}(j_n(\bar{u}_i), \bar{a})=0\ee  
 where $\bar{P}_{i,n}(\bar{x}, \bar{y})$ is over $\kk_0$ and defines the locus of $\la j_n(\bar{u}_i),\bar{a}\ra$ over $\kk_0.$ Since by construction $\bar{a}=\pr_{n,1}(j_n(\bar{u}_i))$, for the regular map $\pr_{n,1}: \Y(n)\to \Y(1)$ over $\Q,$ we can replace (\ref{eqP}) by the equivalent
\be \label{eqP*}\bar{P}^*_{i,n}(j_n(\bar{u}_i))=0, \ee 
where $\bar{P}^*_{i,n}(\bar{x})=\bar{P}_{i,n}(\bar{x}, \pr_{n,1}(\bar{x}))$ defines the locus $Z_{i,n}$ of  $j_n(\bar{u}_i)$ over $\kk_0.$

The systems of 
Zariski closed subsets $Z_{i,n}\subseteq \Y(n)^m$ 
 have the property that
$$\pr_{n,m}(Z_{i,n})=Z_{i,m}, \mbox{ for }m|n.$$
It follows that for some $n_0,$ for any $n,m\ge n_0$
$$\pr_{n,m}\inv(Z_{i,m})=Z_{i,n}, \mbox{ for }m|n,$$
and thus $$\mathrm{locus}(\bar{u}/\kk_0)=j_{n_0}\inv(Z_{1,n_0}),$$
which proves the Lemma for $\kk=\kk_0.$

Note that since $\tilde{\U}_\F$ is invariant under the action of $\Gal(\F/\Q),$ for any $\sigma\in \Gal(\F/\Q),$ 
$$\pr_{n,m}\inv(Z_{i,m}^\sigma)=Z_{i,n}^\sigma, \mbox{ for }m|n,$$

Now let $$Z_{1,n_0}^*=\bigcup_{\sigma\in \Gal(\F/\Q)} Z_{1,n_0}^\sigma. $$
 This is defined over $\Q$ and 
 $$\mathrm{locus}(\bar{u}/\Q)=j_{n_0}\inv(Z_{1,n_0}^*).$$

Now consider $\F_0\subset \F $  algebraically closed and the assumption
\linebreak
$\trd_{\F_0}(\bar{a})>0.$ For this case we refer to the argument  \cite{DH}, section 5.1(b) which show how to deduce the proof of the algebraically  closed case from the case of $\kk_0.$ $\Box$ 

\medskip

{\bf Remark.} In case $j(\bar{u})\in \F_0^k,$ 
$$\mathrm{locus}(\bar{u}/\F_0)=\{ \bar{u}\}.$$

\epk

\bpk \label{k0} {\bf Remark.} Note that  in the first part of the proof  we could set $\kk_0=\Q^\mathrm{ab}(\mathrm{CM}),$ where $\mathrm{CM}$ is the co-ordinates of  all the CM-points on all $\Y(n),$  
Indeed, the key property of such a $\kk_0$ is the one proved in Lemma 5.2 of \cite{DH}: $\kk_0(\bar{a})$ is an abelian extension of $\Q^\mathrm{ab}(\bar{a}).$ This property remains true as we exchange $\mathrm{CM}(1)$ by $\mathrm{CM}$ by the same argument of class field theory.

Moreover, the first statement of Lemma \ref{L6.2} can be extended to the more general:
$$\mathrm{locus}(\bar{u}/\kk)=j_n\inv(Z_{\bar{u},\kk}),$$
where $\kk\subseteq \kk_0,$ any subfield and $Z_{\bar{u},\kk}$ is defined over $\kk.$ The proof of this is the same as the proof for $\kk=\Q.$
\epk

We keep the notation $\F_0$ for an algebraically closed subfield below.

\bpk {\bf Lemma.} 
{\em Let $u_1,\ldots, u_k,v_1,\ldots,v_k\in \U,$ 
$\bar{u}:=\la u_1,\ldots, u_k\ra$ is $\G$-generic, 
$u_i\neq v_i$ $\bar{v}:=\la v_1,\ldots, v_k\ra$ and 
$\bar{v}=\mathsf{g}(\bar{u})$ for  $\mathsf{g}\in \G^k$ (that is
$v_i=g_i\cdot v_i$ for
$g_1,\ldots,g_k\in \G$ respectively, $\mathsf{g}=\la g_1,\ldots,g_k\ra$).

 Then

(i)
$$\mathrm{locus}(\bar{u}\bar{v} /\Q)=
(S\times \U^k)\cap (\mathrm{graph}(\mathsf{g})\cup\mathrm{graph}(\hat{\mathsf{g}}))$$
where $S=\mathrm{locus}(\bar{u}/\Q).$
equivalently
$$\mathrm{locus}(\bar{u}\bar{v} /\Q)=\{ \la x_1,\ldots,x_k,y_1,\ldots,y_k\ra\in \U^{2k}: \vDash S(x) \ \& \ (y=\mathsf{g}(x)\ \vee y=\hat{\mathsf{g}}(x))\}$$

(ii) 
$$\mathrm{locus}(\bar{u}\bar{v} /\F_0)=
(S\times \U^k)\cap \mathrm{graph}(\mathsf{g})$$
where $S=\mathrm{locus}(\bar{u}/\F).$
}

{\bf Proof.} It follows from  the axiom (\ref{B1}) that there are $\mu_1,\ldots,\mu_k\in \hat{\Z}^\times$ such that $$\bar
u\bar{v}\in C_{\gh_1}^{\mu_1}\times \ldots \times C_{\gh_k}^{\mu_k} \subseteq \mathrm{locus}(\bar{u}\bar{v} /\Q)\subseteq C_{\gh_1}\times \ldots\times C_{\gh_k}\mbox{ and } C_{\gh_i}^{\mu_i}=\mathrm{graph}(g_i).$$ 

Taking into account that the $C_{\gh_i,n}^{\mu_i}\subseteq \Y(n)^2$ are algebraic curves defined over $\bar{\Q},$ we get (ii).

By \ref{gg'} $\mathrm{graph}(\gh) \cup \mathrm{graph}(\hat{\gh})$ is the smallest $\Gal_\Q$-invariant subset of $\U^{2k}$ containing  $C_{\gh_1}^{\mu_1}\times \ldots \times C_{\gh_k}^{\mu_k}.$ The statement (i) of Lemma follows. $\Box$

\epk
\bpk\label{L6.3} {\bf Lemma.} {\em Let $\bar{t}=\la t_1,\ldots,t_l\ra\in \U^l$ be the fixed point of
$\mathsf{h}=\la h_1,\ldots,h_l\ra\in \G^l.$  Then
$$\mathrm{locus}(\bar{t} /\F_0)=\{ \bar{t}
\}\mbox{ and }\mathrm{locus}(\bar{t} /\Q)=\{ \bar{t},
\bar{t}^*\}$$
where $\bar{t}^*$
is the unique fixed point of $\hat{\mathsf{h}}:=\la \hat{h}_1,\ldots,\hat{h}_l\ra $ in $\U^l.$}

{\bf Proof.} The first statement is obvious. The second follows from
  \ref{L6.3}, since $\bar{t}^*$ is the only point in $\U^l$ Galois conjugated to $\bar{t}.$ $\Box$

\epk

\bpk\label{L6.5} {\bf Lemma.} {\em Let $u_1,\ldots, u_k, u'_1\ldots,u'_m,v_1,\ldots,v_k, v'_1,v'_l\in \U,$ 
$\bar{u}=\la u_1,\ldots, u_k\ra,$   $\bar{u}'=\la u'_1\ldots,u'_m\ra,$
$\bar{u}\bar{u}'$ is $\G$-generic,
and $S=\mathrm{locus}(\bar{u}\bar{u}').$ Suppose
$u_i\neq v_i,$ for each $1\le i\le k,$ and 
$\bar{v}=\mathsf{g}(\bar{u})$ 
for
$\mathsf{g}\in \G^k.$  Suppose also that $v'_1,\ldots, v'_l$ are special.

Then
$$\mathrm{locus}(\bar{u}\bar{u}'\bar{v}\bar{v}' /\Q)=\{ \bar{x}\bar{x}'\bar{y}\bar{y}'\in \U^{2k+m+l}: \  S(\bar{x}\bar{x}')\ \& \ (\bar{y}=\mathsf{g}(\bar
{x})\ \vee \bar{y}=\hat{\mathsf{g}}(\bar{x}))\ \& \ \bar{y}'\in \{ \bar{t}, \bar{t}^*\}\} $$

$$\mathrm{locus}(\bar{u}\bar{u}'\bar{v}\bar{v}' /\F_0)=\{ \bar{x}\bar{x}'\bar{y}\bar{y}'\in \U^{2k+m+l}: \  S(\bar{x}\bar{x}')\ \& \ (\bar{y}=\mathsf{g}(\bar
{x})\ \& \ \bar{y}'= \bar{t}\} $$
for some $\bar{t}, \bar{t}^*$ as in \ref{L6.3}.
}

{\bf Proof.} Immediate from \ref{L6.2}-\ref{L6.3}. $\Box$
\epk

Call \be
\label{T} T(\bar{x}\bar{x}'\bar{y}\bar{y}'):= S(\bar{x}\bar{x}')\ \& \ (\bar{y}=\mathsf{g}(\bar
{x})\ \vee \bar{y}=\hat{\mathsf{g}}(\bar{x}))\ \& \ \bar{y}'\in \{ \bar{t}, \bar{t}^*\}\ee
{\bf basic predicate} over $\Q.$ 
And
\be \label{TF} T(\bar{x}\bar{x}'\bar{y}\bar{y}'):= S(\bar{x}\bar{x}')\ \& \ \bar{y}=\mathsf{g}(\bar
{x})\ \& \ \bar{y}'= \bar{t}\ee
{\bf basic predicate} over $\F_0,$ where $S,$  $\mathsf{g},$
$\hat{\mathsf{g}},$ $\bar{t}$ and $\bar{t}^*$ are as in \ref{L6.2} -- \ref{L6.5}. 

Clearly, basic predicate define closed analytic subsets.

\bpk \label{6.6} {\bf Proposition.} {\em Let $\kk$ be $\Q$ or an algebraically closed subfield of $\F,$
$T(\bar{x}\bar{x}'\bar{y}\bar{y}')$ be a basic predicate over $\kk,$ $z$ be one of the $M$ variables $x_i, x'_j,y_i$ or $y'_p.$ Let $\pr T\subseteq \U^{M-1}$ be the subset defined by the formula $\exists z\ T(\bar{x}\bar{x}'\bar{y}\bar{y}').$

Then there is a basic analytic set $R$ and a closed subset $R'\subset R,$ both over $\kk,$
$\dim R'< \dim R,$ such that
$$R\setminus R'\subseteq \pr T\subseteq R.$$

}

{\bf Proof.} Let first $\kk:=\Q.$
We consider four possible cases.

(a) $z=y'_n.$ It is immediate from the form of the predicate in (\ref{T})
that $$\exists z\ T(\bar{x}\bar{x}'\bar{y}\bar{y}')\equiv S(\bar{x}\bar{x}')\ \& \ (\bar{y}=\mathsf{g}(\bar
{x})\ \vee \bar{y}=\hat{\mathsf{g}}(\bar{x}))\ \& \ \bar{y}'_-\in \{ \bar{t}_-, \bar{t}^*_-\}$$
where $\bar{y}'_-,  \bar{t}_-$ and $\bar{t}^*_-$ stand for the tuples with omitted $n$-coordinate.

(b) $z=y_i.$ In this case, since  $g_k$ and $\hat{g}_k$ are operations on the whole of $\U,$ 
$$\exists z\ T(\bar{x}\bar{x}'\bar{y}\bar{y}')\equiv  S(\bar{x}\bar{x}')\ \& \ (\bar{y}_- =\mathsf{g}_-(\bar
{x} _-)\ \vee \bar{y}_- =\hat{\mathsf{g}}_-(\bar{x}_-))\ \& \ \bar{y}'\in \{ \bar{t}, \bar{t}^*\}$$
where $\bar{y}_-,$ $\mathsf{g}_-$ and $\bar{x}_-$ stand for the tuples with omitted $i$ - coordinates.

(c) $z=x_i,$ say $i=k.$ This is the same case as (b) if we rearrange the variables in $T$ by taking $x_1,\ldots, x_{k-1},y_k$ to be the set of variables standing for the $\G$-generic $k$-tuple $u_1,\ldots,u_{k-1},v_k$ and replace $\la g_1,\ldots, g_{k-1}, g_k\ra$ by
$\la g_1,\ldots, g_{k-1}, g_k\inv\ra.$ 

(d) $z=x'_j.$ Then 
$$\exists z\ T(\bar{x}\bar{x}'\bar{y}\bar{y}')\equiv (\exists x'_j\ S(\bar{x}\bar{x}'))\ \& \ (\bar{y}=\mathsf{g}(\bar
{x})\ \vee \bar{y}=\hat{\mathsf{g}}(\bar{x}))\ \& \ \bar{y}'\in \{ \bar{t}, \bar{t}^*\}.$$

Let $\pr S$ be the subset of $\U^{k+m-1}$ defined by $\exists x'_j\ S(\bar{x}\bar{x}')$ and
$L$ the subset of $\U^{2k+l}$ defined by $(\bar{y}=\mathsf{g}(\bar
{x})\ \vee \bar{y}=\hat{\mathsf{g}}(\bar{x}))\ \& \ \bar{y}'\in \{ \bar{t}, \bar{t}^*\}.$

Since $S=j_n\inv(W)$ for some  Zariski closed subset $W\subseteq \X(n)^{k+m}, $ 
we have $\pr S= j_n\inv(\pr W),$ where on the right of the equation we consider the projection along $x'_j$ and on the left the projection along the coordinate corresponding in the image. By standard facts on Zariski topology 
$\pr W= V\setminus R,$
for $V$ Zariski closed and $U$ a boolean combination of Zariski closed,
$\dim U<\dim V.$  Let $\bar{U}$ be the Zariski closure of $U,$ so
the Zariski open set $V\setminus \bar{U}$ is a subset of $\pr W$ and $\pr W\subset V.$

Then $\pr S=j_n\inv(\pr W)=j_n\inv(V)\setminus j_n\inv(U)$ and so
$$R\setminus R'\subseteq \pr S\subseteq R, \mbox{ for } R=j_n\inv(V),\ R'=j_n\inv(\bar{U}).$$

It follows that, 
$$((R\setminus R')\times \U^{k+l})\cap (\U^m\times L) \subseteq 
\exists z\ T(\bar{x}\bar{x}'\bar{y}\bar{y}')\subseteq (R\times \U^{k+l})\cap (\U^{m-1}\times L)$$
which proves the Proposition for $\kk=\Q.$

For $\kk=\F_0$ algebraically closed, use the same arguments (a)-(d) combined with $\F_0$-versions of \ref{L6.2}-\ref{L6.5}.

$\Box$

\medskip

Let $T\subseteq \U^M$ be a basic predicate over $\kk,$ that is a predicate of the form
 (\ref{eqP}) or (\ref{eqP*}).
 $$\dim T:\dim S:=\dim Z$$ 
 where $Z\subseteq \Y(n)^{k+m}$ is the Zariski closed subset over $\kk$ such that, according to  \ref{L6.2}, 
 $S=j_n\inv(Z).$

\epk 
\bpk\label{P6.7} {\bf Proposition.} {\em Suppose $\F$ is of infinite transcendence degree over $\kk,$ for $\kk=\Q$ or $\kk=\F_0,$ an algebraically closed subfield of $\F.$ Then 
$\U_\F$ in the language of basic predicates over $\kk$ is $\omega$-homogeneous over $\kk$:

for any $\bar{u},\bar{u}'\in \U^m$ and $v\in \U$ such that 
$ \mathrm{locus}(\bar{u}/\kk)=\mathrm{locus}(\bar{u}'/\kk)$
there is $v'\in \U$ such that
$\mathrm{locus}(\bar{u}v/\kk)=\mathrm{locus}(\bar{u}'v'/\kk).$

}

{\bf Proof.} Let $T=\mathrm{locus}(\bar{u}v/\kk)$ and  $R=\mathrm{locus}(\bar{u}/\kk).$ Then by \ref{6.6} there is a closed $R'$ of smaller dimension such that 
$$R\setminus R'\subseteq \pr T\subseteq R.$$
Clearly, $\bar{u}'\in R\setminus R'$ and hence $\bar{u}'\in \pr T,$ which means that 
the fibre
$$T(\bar{u}',\U):=\{ w\in \U: \   \bar{u}'w\in T\}\neq \emptyset.$$

Note that $T(\bar{u}',\U)$ is a closed subset of $\U$ defined over the  field $\kk'$ generated by coordinates of $j_n(\bar{u}'),$ all $n\in \N.$

Consider the only two possible cases, $\dim T=\dim R$ and $\dim T=\dim R+1.$ 

In the first case pick up any $v'\in T(\bar{u}',\U).$ Then
$\mathrm{locus}(\bar{u}'v'/\Q)=T$ since $T$ is irreducible over $\Q.$
   
In the second case $T(\bar{u}',\U)=\U.$ Pick up any $v'\in \U$ generic over $\kk',$ that is such that $j_1(v')\notin \bar{\kk}'.$
 Clearly, the transcendence degree of $\kk'$ over $\Q$ is at most the length of $\bar{u},$ and so such $v'$ exists. Now
 $\dim \mathrm{locus}(\bar{u}'v'/\Q)=\dim R+1$ and so again by irreducibility of $T$ we have the equality $\mathrm{locus}(\bar{u}'v'/\Q)=T.$ $\Box$
 

\epk
\bpk For a subset $W\subset \U$ and a point $u\in \U$ we define
$$u\in \cl(W)\Leftrightarrow \exists \bar{w}\subset_\mathrm{finite} W:
\dim \mathrm{locus}(\bar{w}u)=\dim \mathrm{locus}(\bar{w}).$$
And define the closure of $W$
$$\cl(W)=\{ u\in \U: \ u\in \cl(W)\}.$$

We consider the covering sort $\U$ in $\U_\F$ together with basic relations over $\Q$ as a structure. 

Recall (see \cite{K5}) that one calls $(\U,\cl)$ a  {\bf quasiminimal pregeometry structure} if the following holds:

QM1. The pregeometry is determined by the language. That is, if $\tp(v\bar{w}) = \tp(v'\bar{w}')$
 then $v\in  \cl(\bar{w})$ if and only if $v'\in  \cl(\bar{w}').$

QM2. $\U$ is infinite-dimensional with respect to $\cl.$

QM3. (Countable closure property) If $W \subset  \U$ is finite then $\cl(W)$ is countable.

QM4. (Uniqueness of the generic type) Suppose that $W, W'\subseteq \U$ are countable
subsets, $\cl(W)=W,$ $\cl(W')=W'$ and $W,W'$ enumerated so that $\tp(W) = \tp(W').$ 

If $v \in \U \setminus W$
and $v'\in  \U \setminus W'$ then $\tp(W v) = \tp(W' v')$
 (with respect to the same
enumerations for $W$ and $W'$).

QM5. ($\aleph_0$-homogeneity over closed sets and the empty set)
Let $W, W' \subseteq \U$ be countable closed subsets or empty, enumerated such
that $\tp(W) = \tp(W'),$ and let $\bar{w}, \bar{w}'$
 be finite tuples from $\U$ such that
$\tp(W\bar{w}) = \tp(W'\bar{w}')$, and let $v \in \cl(W\bar{w}).$ Then there is $v'\in \U$ such
that $\tp( \bar{w} vW) = \tp(\bar{w}' v'W').$

\epk

\bpk \label{MThm} {\bf Theorem.} {\em For any model $\U_\F$ of $\Sigma,$
with $\F$ algebraically closed of infinite transcendence degree, the structure $(\U,\cl)$ is a quasiminimal pregeometry. 

There is an existential $L_{\omega_1,\omega}$-sentence $\Theta_{\U}$ such that $\Sigma  \& \Theta_{\U}$
 defines a categorical AEC containing $\U.$  
 }

{\bf Proof.} We strat with the proof of the firs statement of the theorem by checking conditions QM1--QM5.

QM1 is by definition. QM2 follows fom the assumption on $\F.$

QM3 is due to the fact, implied by the definition, that $v\in \cl(W)$ if and only if $j(v)\in \acl(j(W)),$ taking into account that $j\inv(j(v))$ is countable.

To tackle QM4 and QM5 note first that $\tp(\bar{u})$ is determined by
$T=\mathrm{locus}(\bar{u}/\Q).$ More presicely, the quantifier-free part of the type is given by the basic predicate $T$ together with negations of all the basic predicates $R$ of the same arity such that $\dim R<\dim T.$  Now we claim that any type is equivalent to a quntifier-free one. Indeed, by homogeneity proved in \ref{P6.7}, in a countable elementary substructures $\U^0\prec \U$ the set defined by
 $\tp(\bar{u})$ is equal to the set defined by the respective quantifier-free type. Hence the same holds in $\U.$
 
 Now note that the condition $v\notin W=\cl(W)$ in QM4 is equivalent to
the condition that $j_n(v)$ is generic in $\Y(n)$ over $j_n(W)),$ for all (equivalently, for some) $n\in \N.$ Since $\Y(n)$ is absolutely irreducible, the condition determines the complete field-theretic type of $j_n(v)$ over $W$ and hence the complete  quantifier-free type of $v$ over $W$, equivalently, the full type of $v$ over $W.$ QM4 follows.

QM5 is a direct consequence of \ref{P6.7}.
This completes the proof of the first statement.

Now we construct the $L_{\omega_1,\omega}$-sentence $\Theta_\U.$ 

For each tuple $\mathsf{h}=\la h_1,\ldots h_l\ra\in \E^l$ and  the respective tuple $\bar{t}_h$ of fixed element $\la t_1,\ldots,t_l\ra$ of $\mathsf{h},$ for each $n$
 consider the minimal Zariski closed subset $Z_{h,n}\subset \Y(n)^l$ over $\Q$ such that $j_n(\bar{t}_h)\in Z_{h,n}.$ These depend on the model $\U_\F$ of $\Sigma.$
  
  Now set  $$\Phi_h=\bigwedge_{n\in \N} j_n(\bar{t}_h)\in Z_{h,n}.$$ This can be seen as a quantifier-free   $L_{\omega_1,\omega}$-formula     
with variables $\bar{t}_h.$ 

 $\Theta_\U$ will be the 
 $L_{\omega_1,\omega}$-formula $\Theta$ constructed in \ref{listG}  
  stating that
for any $\mathsf{h}\in \E^l$ there exists $\bar{t}_h=\la t_1,\ldots,t_l\ra,$ a tuple of fixed points of $\mathsf{h}$ such that $\Phi_i(\bar{t}_h).$ 

The rest of the proof is split into Lemmas and Claims below.

 Now let $\F$ be an uncountable algebraically closed field of characteristic zero consider $\U_\F$ and $\U'_F$ models  of $\Sigma \ \& \ \Theta.$ We aim to prove that $\U_\F\cong\U'_F.$

\medskip

Let $\G$ and $\G'$ be the realisations of $\GL_2^+(\Q)$ in $\U_\F$ and $\U'_\F$ respectively. By construction, the statement in the formula $\Theta$ implies that some group-isomorphism $\mathsf{i}_\G: \G\to \G'$ can be uniquely extended to the map  $\mathsf{i}_\mathrm{CM}: \U(\mathrm{CM})\to \U'(\mathrm{CM})$  which acts on the fixed points 
$u_g\mapsto u_{g'},$ if $g'=\mathsf{i}_\G(g)$ and takes
$j_n(u_g)$ to  $j_n(u_{g'})$ so that
 polynomial equations in the co-ordinates of $j_n(u_g)$ over $\Q$ are preserved. In other words we have a Galois automorphism   $$\mathsf{i}_{\kk}: \Q(\mathrm{CM})\to \Q(\mathrm{CM}), \mbox{ where }\kk:=\Q(\mathrm{CM}), $$
which agrees with $\mathsf{i}_\G,$ that is together the pair $(\mathsf{i}_\G,
 \mathsf{i}_{\kk})$ is a partial isomorphism $\mathsf{i}_0:\U_\F\to \U'_\F.$
 
Let $\F_0\subseteq \F$ be a countable algebraically closed subfield. 

We consider submodels $\U_{\F_0}\subseteq \U_\F$ and  $\U'_{\F_0}\subseteq \U'_\F.$  
 
Claim 1. There are subsets $V\subset  \U_{\F_0}$
and $V'\subset \U'_{\F_0}$
and a partial isomorphism $\mathsf{i}_V: V\to V'$ extending $\mathsf{i}_0$
such that
 $V$ is a maximal $\G$-generic subset of $\U_{\F_0}$  and 
$V'$ is maximal $\G$-generic subset of $\U'_{\F_0}.$

 Proof. Since $\U_{\F_0}$ and $\U'_{\F_0}$ are countable we can enumerate the sets and use the back-and-forth procedure in constructing
 $V=\{ v_m:m\in \N\},$ $V'=\{ v'_m:m\in \N\},$ and $\mathsf{i}_V.$
 
 Suppose $\bar{v}=\la v_1,\ldots,v_m\ra$   and $\bar{v}'=\la v'_1,\ldots,v'_m\ra$ be $\G$-generic  and satisfy
 $$\mathrm{locus}(\bar{v}'/\kk)=\mathsf{i}_\kk(\mathrm{locus}(\bar{v}/\kk)).$$
 which means that $j_n(\bar{v})\in Z_n,$ for some variety $Z_n\subset \Y(n)^m$ over $\kk,$ if and only if $j_n(\bar{v}')\in Z'_n$ for the variety $Z'_n=\mathsf{i}_\kk(Z_n),$ for all $n.$
 
For the back-and-forth construction, it suffices to prove: 

Let $w$ be the first element in $\U_{\F_0}$ such that 
$\bar{v}w$  is $\G$-generic. There exists $w'\in \U'_{\F_0}$ such that
$$\mathrm{locus}(\bar{v}'w'/\kk)=\mathsf{i}_\kk(\mathrm{locus}(\bar{v}w/\kk)).$$
In order to establish $w'$ consider $\{ W_n\subset \Y(n)^{m+1}: n\in \N\}$ be the family of varieties over $\kk$ which determine the locus of $\bar{v}w.$ By \ref{L6.2} together with \ref{k0} there is $n_0$ such that the locus is actually determined by any one of the $W_n,$ for $n\ge n_0.$ By algebraic geometry
the projection of $\pr W_n$ on the first $m$ coordinates contains a Zariski open subset of  $Z_n,$ hence contains $j_n(\bar{v}).$
 Since $Z_n$ is Galois conjugated to $Z'_n,$ there exists respective $W'_n=\mathsf{i}_\kk(W_n)$
together with the whole family  $\{ W'_n\subset \Y(n)^{m+1}: n\in \N\}$
such that $j_n\inv(W'_n)=j_{n_0}\inv(W'_{n_0}),$ for $n\ge n_0,$
and  $\pr W'_n$ contains an open subset of $Z'_n.$
Thus $j_n(\bar{v}')=\bar{z}'\in \pr W'_n$ that is there exists $t\in \Y(n),$ such that
$\bar{z}t\in W'_n,$ generic over $\kk.$ Pick up $w'\in j_n\inv(t),$ and thus
$j_n(\bar{v}'w')\in W'_n.$ So $w'$ is as required.
Claim 1 proved.
 
 \medskip
 
Claim 2.   There is a unique extension $\mathsf{i}$ of $\mathsf{i}_V$ to
an isomorphism
$$\mathsf{i}: \U_{\F_0}\to \U'_{\F_0}.$$

Proof. By definition $$\U_{\F_0}=\U(\mathrm{CM})\ \dot{\cup}\ \G\cdot V\mbox{ and }\U'_{\F_0}=\U'(\mathrm{CM})\ \dot{\cup}\ \G'\cdot V',$$
where $$\G\cdot V=\{ \G\cdot v: v\in V\},\ \ \G'\cdot V'=\{ \G'\cdot v': v'\in V'\}$$
the disjoint unions of Hecke orbits of non-CM points. Any $\mathsf{i}$ must coincide with $\mathsf{i}_\mathrm{CM}$ on $\U(\mathrm{CM}),$ and on the remaining points define
 $$\mathsf{i}: g\cdot v\mapsto \mathsf{i}_\G(g)\cdot \mathsf{i}_V(v),$$
 $$\mathsf{i}: j_n(g\cdot v)\mapsto j_n(\mathsf{i}_\G(g)\cdot \mathsf{i}_V(v)).$$
This is as required.   

\medskip

Finally we are ready to prove  $\U_\F\cong \U'_\F.$ 

Let $\F_\omega\subset \F$ be an algebraically closed subfield of an infinite countable transcendence degree. 
We follow \cite{K5}
and consider the  class $\mathcal{K}(M)$ (see Theorem 2.3 therein and the definition before it) for $M\cong\U_{\F_\omega}\cong \U'_{\F_\omega}.$ 
The class of countable substructures of $\U_\F$ coincides, up to isomorphism, to the   class of countable substructures of $\U'_\F,$ by Claim 2.
Thus $\U_\F$ and $\U'_\F$ belong to $\mathcal{K}(\U_{\F_\omega}),$ by definition.

Since
$\U_{\F_\omega}$ is a quasiminimal pregeometry structure,  $\mathcal{K}(\U_{\F_\omega})$ is categorical in uncountable cardinalities. The statement follows.
$\Box$ 
\epk
\pagebreak

\end{document}